\documentclass[11pt]{article}
\title{On Ikehara type Tauberian theorems with $O(x^\gamma)$ remainders} 
\author{Michael M\"uger \\ Institute for Mathematics, Astrophysics and Particle Physics \\
Radboud University, Nijmegen, The Netherlands}

\usepackage{amssymb,amsmath, theorem,latexsym,ifthen,comment}
\usepackage{mathrsfs}
\usepackage{epic}
\usepackage{eepicemu}
\usepackage{epsfig}
\newboolean{noremarks}
\newlength{\dinwidth}
\newlength{\dinmargin}
\def\1#1{{\bf #1}}
\def\2#1{{\cal #1}}
\def\3#1{{\sl #1}}
\def\4#1{{\tt #1}}
\def\5#1{{\sf #1}}
\def\6#1{{\mathfrak #1}}
\def\7#1{{\mathbb #1}}
\def\8#1{{\mathscr #1}}

\newcommand{\be}{\begin{equation}}
\newcommand{\ee}{\end{equation}}
\newcommand{\ba}{\begin{array}}
\newcommand{\ea}{\end{array}}
\newcommand{\bea}{\begin{eqnarray}}
\newcommand{\eea}{\end{eqnarray}}
 \newcommand{\bean}{\begin{eqnarray*}}
\newcommand{\eean}{\end{eqnarray*}}
\newcommand{\nn}{\nonumber}

\newcommand{\ve}{\varepsilon}

\newcommand{\impl}{\Rightarrow}
\newcommand{\rarr}{\rightarrow}

\renewcommand{\Re}{\mathrm{Re}}

\newcommand{\qed}{\hfill$\blacksquare$\\}

\def\endexem{\hfill{$\Box$}\medskip}

\theoremstyle{change}
\theoremheaderfont{\scshape}
\newtheorem{defin}{Definition}[section]
\newtheorem{defprop}{Definition/Proposition}[section]
\newtheorem{lemma}[defin]{Lemma}
\newtheorem{prop}[defin]{Proposition}
\newtheorem{theorem}[defin]{Theorem}
\newtheorem{coro}[defin]{Corollary}
\newtheorem{conj}[defin]{Conjecture}
\newtheorem{wish}[defin]{Desideratum}
\newtheorem{question}[defin]{Question}
\theorembodyfont{\rmfamily}
\newtheorem{remark}[defin]{Remark}
\newtheorem{exercise}[defin]{Exercise}

\newcommand{\bdefin}{\begin{defin}}
\newcommand{\blemma}{\begin{lemma}}
\newcommand{\bprop}{\begin{prop}}
\newcommand{\btheor}{\begin{theorem}}
\newcommand{\bcoro}{\begin{coro}}
\newcommand{\bdefprop}{\begin{defprop}}
\newcommand{\bexer}{\begin{exercise}}
\newcommand{\edefprop}{\end{defprop}}
\newcommand{\edefin}{\end{defin}}
\newcommand{\elemma}{\end{lemma}}
\newcommand{\eprop}{\end{prop}}
\newcommand{\etheor}{\end{theorem}}
\newcommand{\ecoro}{\end{coro}}
\newcommand{\eexer}{\end{exercise}}
\newcommand{\bconj}{\begin{conj}}
\newcommand{\econj}{\end{conj}}
\newcommand{\bwish}{\begin{wish}}
\newcommand{\ewish}{\end{wish}}
\newcommand{\bquestion}{\begin{question}}
\newcommand{\equestion}{\end{question}}
\newcommand{\brem}{\begin{remark}}
\newcommand{\erem}{\endexem\end{remark}}

\newcommand{\prf}{{\noindent\it Proof. }}
\def\mobj#1{\raise .4\unitlens\hbox{\put(0,0){$#1$}}}
\def\mychi{\raise 2pt\hbox{$\chi$}}

\begin{document}

\maketitle

\abstract{Motivated by analytic number theory, we explore remainder versions of Ikehara's Tauberian 
theorem yielding power law remainder terms. More precisely, for $f:[1,\infty)\rarr\7R$ non-negative
and non-decreasing we prove $f(x)-x=O(x^\gamma)$ with $\gamma<1$ under certain assumptions on $f$. We
state a conjecture concerning the weakest natural assumptions and show that we cannot hope for more.}

\section{Motivation and results}
The following was proven in 1931 by Wiener's student Ikehara \cite{ikehara}. (For a much better proof see
\cite{bochner} or \cite[Section 3.5]{CQ}.)

\btheor [Ikehara, 1931] \label{theor-ikehara}
Let $f:[1,\infty)\rarr\7R$ be non-negative and non-decreasing. Assume that 
\[ F(s)=\int_1^\infty f(x)x^{-s}\frac{dx}{x} \]
converges for $s>1$ (thus $F$ is holomorphic on $\{\Re\, s>1\}$). Assume that $F(s)-\frac{A}{s-1}$ has
a continuous extension to the closed half-plane $\{ \Re\, s\ge 1\}$. Then $f(x)=Ax+o(x)$.
\etheor

This result gives rise to what still is the simplest proof of the prime number theorem
$\pi(x)\sim\frac{x}{\log x}$. In most approaches to giving more precise estimates of
$\pi(x)-\mathrm{Li}(x)$, or rather of $\psi(x)-x$, Tauberian theorems have not played a major
r\^ole. One exception is provided by \cite{kienast, cizek}, where remainder terms of the form
$\frac{x}{\log^kx}$ are proven under somewhat stronger assumptions than in Theorem
\ref{theor-ikehara}, which are then used to give the simplest known proofs of
$\psi(x)-x=O(\frac{x}{\log^kx})\ \forall k\in\7N$, invoking only properties of $\zeta(s)$ for
$\Re\,s\ge 1$. 
More general results on remainder estimates in Ikehara's theorem are found \cite[\S 7.5]{tenenbaum},
\cite{RR}, but the Tauberian conditions considered here are different: 
%We want to prove bounds on $f(x)-Ax$ from existence of holomorphic continuations of
%$F(s)-\frac{A}{s-1}$ to halfplanes $\{\Re\, s>\alpha\}$.  

\bquestion \label{conj}
Assume that $f:[1,\infty)\rarr\7R$ is non-negative and non-decreasing, the integral
$F(s)=\int_1^\infty f(x)x^{-s-1}dx$ converges for $s>1$ (thus for $\Re\,s>1$), and
$F(s)-\frac{A}{s-1}$ has a holomorphic extension to the half-plane $\{ \Re\, s>\alpha\}$, where
$\alpha\in(0,1)$. Does this imply $f(x)=Ax+O(x^{\lambda+\ve})$ for some $\lambda<1$?
\equestion

\brem 1. Ikehara's theorem shows that there is a unique $A$ such that $f(x)-Ax=o(x)$.

2. It is trivial that if $g:[1,\infty)\rarr\7R$ is measurable and $g(x)=O(x^\gamma)$, then
$G(s)=\int_1^\infty g(x)x^{-s-1}dx$ is convergent for $\Re\,s>\gamma$ and defines a holomorphic
function on this domain. The above question is equivalent to asking to which extent this can be
inverted under the additional assumption that $x\mapsto g(x)+Ax$ is non-decreasing for some $A$
(keeping in mind that the domain of holomorphicity of $G$ can be larger then the domain of
convergence of the integral).
 
3. If the answer to the question was positive with $\lambda=\alpha$ in case $A>0$, it would provide
a very simple deduction of $\psi(x)-x=O(x^{\alpha+\ve})$ from 
$\Re\,s>\alpha\impl\zeta(s)\ne 0$. However, we will prove see that this is not the case.  
\erem

But there is a weaker positive answer to the question. To wit, we will prove the following:

\btheor \label{theor-main}
Let $f:[1,\infty)\rarr\7R$ be non-negative and non-decreasing. Assume that 
$F(s)=\int_1^\infty f(x)x^{-s-1}dx$ converges for $s>1$ and that $F(s)-\frac{A}{s-1}$ has a
holomorphic continuation to $\{\Re\,s>\alpha\}$, where $\alpha\in(0,1)$. 
Then $f(x)=Ax+O(x^{\gamma+\ve})$, where
\begin{itemize}
\item[(I)] $\gamma=\alpha$ if $A=0$, 
\item[(II)] $\gamma=\frac{\alpha+1}{2}$ if $A>0$ and $f(x)-Ax$ is of fixed sign for $x\ge x_0$ for
some $x_0$.
\end{itemize}
These exponents are optimal under the given assumptions.
\etheor

These results will be proven in Section \ref{sec2}. In Section \ref{sec3} we show that the statement
in Case II is false without the sign condition. We also conjecture that $\gamma=\frac{\alpha+2}{3}$
always works and provide some evidence.  

\vspace{.5cm}

\noindent{\it Acknowledgment.}  
I thank G.\ Tenenbaum for alerting me \cite{t} of an egregious error in the first version of this
note and for pointing out that Proposition \ref{prop-Omega} is about convergence, not absolute
convergence, which led to a simpler proof.

%%%%%%%%%%%%%%%%%%

\section{Proofs}\label{sec2}
Case (I) in Theorem \ref{theor-main} is fairly trivial and surely well-known. We only include the
proof as a preparation for the following one. 

\blemma \label{lem-f} Let $\alpha>0$ and $f:[1,\infty)\rarr\7R$ be non-negative and
non-decreasing. If $\int_1^\infty f(x)x^{-s-1}dx$ converges for all $s>\alpha$ then
$f(x)=O(x^{\alpha+\ve})$.  
\elemma

\prf Assume $f(x)=\Omega(x^\gamma)$ with $\gamma>\alpha$. Then there are $C>0$ and arbitrarily large
$z$ such that $f(z)\ge Cz^\gamma$. For such a $z$, we have $f(x)\ge f(z)\ge Cz^\gamma$ whenever
$x\ge z$. Taking $s=\gamma$, we have $\int_z^{2z} f(x)x^{-s-1}dx\ge C\int_z^{2z}x^{\gamma-s-1}dx=C\log 2$.
This contradicts the assumed convergence of $\int_1^\infty f(x)x^{-s-1}dx$
(since the latter means that for every $\ve>0$ there is a $T$ such that $T\le x_1\le x_2$ implies
$|\int_{x_1}^{x_2}f(x)x^{-s-1}dx|<\ve$). This contradiction proves that $f(x)=o(x^\gamma)$ for
all $\gamma>\alpha$, which is equivalent to the assertion.
\qed

\brem With $f(x)=x^\alpha\log x$, we have convergence of $\int_1^\infty f(x)x^{-s-1}dx$ for all 
$s>\alpha$, but $f(x)=O(x^\gamma)$ holds if and only if $\gamma>\alpha$. This proves optimality of
the result of the lemma.
\erem

\bprop \label{prop-Omega}
Let $f:[1,\infty)\rarr\7R$ be non-negative and non-decreasing. Then
\begin{itemize}
\item[(i)] If $A>0$ and $\gamma\in(0,1)$ are such that $f(x)-Ax=\Omega(x^\gamma)$ then 
$\int_1^\infty (f(x)-Ax)x^{-s-1}dx$ diverges whenever $s\le 2\gamma-1$.
\item[(ii)] Let $\sigma_c$ be the abscissa of conditional convergence of 
$\int_1^\infty(f(x)-Ax)x^{-s-1}dx$. Then $f(x)-Ax=O(x^{\gamma+\ve})$,
where $\gamma=\frac{\sigma_c+1}{2}$.
\end{itemize}
\eprop

\prf (i) The assumption $f(x)-Ax=\Omega(x^\gamma)$ means that there are $C>0$ and arbitrarily large
$x$ such that $|f(x)-Ax|\ge Cx^\gamma$. Assume $x_1$ is such that $f(x_1)-Ax_1\ge Cx_1^\gamma$. 
Since $f$ is non-decreasing, we have $f(x)\ge f(x_1)\ge Ax_1+Cx_1^\gamma$ for all $x\ge x_1$. Put  
$x_2=x_1+\frac{Cx_1^\gamma}{2A}$. Then for all $x\in[x_1,x_2]$ we have 
\[ f(x)-Ax\ge Ax_1+Cx_1^\gamma-Ax_2= Cx_1^\gamma-A(x_2-x_1)= Cx_1^\gamma-A\frac{Cx_1^\gamma}{2A}
   =\frac{Cx_1^\gamma}{2}.\]
Thus for $s\ge -1$ we have
\[ \int_{x_1}^{x_2}\frac{f(x)-Ax}{x^{s+1}}dx \ge 
   \frac{Cx_1^\gamma}{2A} \frac{Cx_1^\gamma}{2}\frac{1}{x_2^{s+1}}
   =\frac{C^2x_1^{2\gamma}}{4A(x_1+\frac{Cx_1^\gamma}{2A})^{s+1}}
   =\frac{C^2}{4A}\frac{x_1^{2\gamma-s-1}}{(1+\frac{C}{2Ax_1^{1-\gamma}})^{s+1}}. \]
Now assume $x_2$ is such that $f(x_2)-Ax_2\le -Cx_2^\gamma$. Since $f$ is non-decreasing, this
implies that $f(x)\le Ax_2-Cx_2^\gamma$ for all $x\le x_2$. Define $x_1=x_2-\frac{Cx_2^\gamma}{2A}$. 
By a reasoning similar to the one above we have $f(x)-Ax\le -\frac{Cx_2^\gamma}{2}$ for all 
$x\in[x_1,x_2]$, implying 
\[ \int_{x_1}^{x_2}\frac{f(x)-Ax}{x^{s+1}}dx \le
   -\frac{C^2}{4A}\frac{x_2^{2\gamma-s-1}}{(1-\frac{C}{2Ax_2^{1-\gamma}})^{s+1}}. \]
If $s\le 2\gamma-1$ then $2\gamma-s-1\ge 0$. Then the above computations and
$f(x)-Ax=\Omega(x^\gamma)$ imply that for every $T\ge 1$ we can find an interval
$[x_1,x_2]\subseteq [T,\infty)$ such that 
\[ \left|\int_{x_1}^{x_2}\frac{f(x)-Ax}{x^{s+1}}dx \right| \ge 
 \frac{C^2}{4A \cdot 2} x_1^{2\gamma-s-1} \ge \frac{C^2}{4A \cdot 2}. \]
But this clearly implies that $\int_1^\infty(f(x)-Ax)x^{-s-1}dx$ diverges.

(ii) Follows from the combination of (i) with convergence of the integral for each $s>\sigma_c$.
\qed

Note that the above does not assume $f(x)-Ax$ to be ultimately of constant sign.

The following is a version of the Phragm\'en-Landau theorem on Dirichlet series with positive 
coefficients, cf.\ e.g.\ \cite[Sec. II.1, Theorem 1.9]{tenenbaum}.

\bprop \label{prop-conv} Assume that $g:[1,\infty)\rarr\7R$ is non-negative and measurable, that
\be G(s)=\int_1^\infty g(x)x^{-s}\frac{dx}{x} \label{G}\ee
converges for $\Re\,s>1$ and that the function $G$ has a holomorphic extension to the 
half-plane $\{ \Re\, s>\alpha\}$, where $\alpha\in(0,1)$. Then the integral in (\ref{G}) converges
to $G$ whenever $\Re\,s>\alpha$.
\eprop

%\begin{comment}
\prf Let $\alpha<t<1<s$. Since $G$ is holomorphic at $s$, it has a power series expansion
\be G(z)=\sum_{n=0}^\infty \frac{(z-s)^n\,G^{(n)}(s)}{n!}.\label{G1}\ee
Since $\int_1^\infty g(x)x^{-s}\frac{dx}{x}<\infty$ for all $s>1$ and 
$(\log x)^k=O(x^\ve)$ for any $k,\ve>0$, we also have 
$\int_1^\infty|g(x)|(\log x)^nx^{-s}\frac{dx}{x}<\infty$ for all $s>1, n\in\7N$. With 
$\frac{d^n}{ds^n}x^{-s}=(-\log x)^nx^{-s}$ and Lebesgue's dominated convergence theorem we can
differentiate $G(s)$ under the integral sign and obtain
\[ G^{(n)}(s)=\int_1^\infty g(x)(-\log x)^n x^{-s}\frac{dx}{x}. \]
Since $G$ is holomorphic on the half-plane $\{\Re\,z>\alpha\}$, the domain of convergence of (\ref{G1})
includes $t$. Thus
\[ G(t)=\sum_{n=0}^\infty \frac{(t-s)^n}{n!}\int_1^\infty g(x)(-\log x)^n x^{-s}\frac{dx}{x}
  =\sum_{n=0}^\infty \int_1^\infty \frac{(s-t)^n}{n!}g(x)(\log x)^n x^{-s}\frac{dx}{x}. \]
Since the integrand is non-negative and the double integral converges by our assumptions, by
Fubini-Tonelli we may reverse the order of summation and integration: 
\bean G(t) &=& \int_1^\infty \left(\sum_{n=0}^\infty \frac{(s-t)^n}{n!} (\log x)^n\right)g(x)x^{-s}\frac{dx}{x}
   =\int_1^\infty e^{(s-t)\log x} g(x)x^{-s}\frac{dx}{x} \\
  &=& \int_1^\infty x^{s-t} g(x)x^{-s}\frac{dx}{x} =\int_1^\infty g(x)x^{-t}\frac{dx}{x},
\eean
where the rightmost integral converges. 
\qed
%\end{comment}

\noindent{\it Proof of Theorem \ref{theor-main}, Case (II).}
Assume that $g(x):=f(x)-Ax$ has constant sign for $x\ge x_0$. Since
$\int_1^{x_0}(f(x)-Ax)x^{-s-1}dx$ converges for all $s\in\7C$ and defines an entire function, we may
replace the lower integration bound $1$ by $x_0$ in the argument that follows, so that $g$ has
constant sign on the domain of integration.
It is clear that $G(s)=\int_{x_0}^\infty g(x)x^{-s-1}dx$ converges for all $s>1$, and the function
$G$ by assumption continues holomorphically to $\{\Re\,s>\alpha\}$. Thus Proposition \ref{prop-conv}
(which of course also holds for non-positive functions) implies that the integral converges to $G$
for $s>\alpha$. Now the claim follows from Proposition \ref{prop-Omega}. 
\qed

That the statements of Proposition \ref{prop-Omega} and Case II of the theorem are optimal follows
from the following example: 

\bprop \label{prop1} For every $\gamma\in(0,1)$ there exists a function $f:[1,\infty)\rarr\7R$ such
that
\begin{itemize}
\item $f$ is non-decreasing and $f(x)\ge x\ \forall x\ge 1$,
\item $f(x)-x$ is $O(x^\gamma)$ and $\Omega(x^\gamma)$ (thus not $O(x^{\gamma'})$ for any
$\gamma'<\gamma$). 
\item $F(s)=\int_1^\infty f(x)x^{-s-1}dx$ converges for $s>1$,
\item $G(s)=\int_1^\infty (f(x)-x)x^{-s-1}dx$ converges if and only if $s>2\gamma-1$,
\item $G(s)=F(s)-\frac{1}{s-1}$ is holomorphic on  $\{\Re\,s>2\gamma-1\}$ and has a singularity at
$2\gamma-1$. 
\end{itemize}
\eprop

\prf Let $\{x_n\}, \{h_n\}$ be sequences satisfying $x_1\ge 1$ and 
$x_{n+1}\ge x_n+h_n\ \forall n$. Define $f:[1,\infty)\rarr\7R$ by
\[ f(x)=\left\{ \begin{array}{ll} x_i+h_i & \mathrm{if}\ x\in[x_i,x_i+h_i] \ \mathrm{for\ some}\ i\\ x &
    \mathrm{otherwise} \end{array}\right. \]
It is obvious that $f$ is non-negative, non-decreasing and satisfies $f(x)-x\ge 0\ \forall x$. With
$h_n=x_n^\gamma$ it is immediate that $f(x)-x$ is $O(x^\gamma)$ and $\Omega(x^\gamma)$ (since
$f(x_i)=x_i+x_i^\gamma$ for all $i$ and $x_i\rarr\infty$).
In view of $0\le f(x)-x\le h_i=x_i^\gamma$ for $x\in[x_i,x_i+h_i]$, we have
\[ G(s)=\int_1^\infty (f(x)-x)x^{-s-1}dx\le\sum_{i=1}^\infty x_i^{2\gamma_i} x_i^{-s-1}. \]
Taking $x_i=2^i$, the r.h.s.\ becomes $\sum_{i=1}^\infty 2^{(2\gamma-s-1)i}$, which converges
whenever $2\gamma-s-1<0$, or $s>2\gamma-1$. Thus for the integral defining $G$ we have 
$\sigma_c=\sigma_a\le 2\gamma-1$, and $G$ is holomorphic on $\{\Re\,s>2\gamma-1\}$.
Proposition \ref{prop-Omega} gives $\sigma_c\ge 2\gamma-1$ and Proposition \ref{prop-conv} implies
that $G$ has a singularity at $\sigma_c$ (which a more careful computation shows to be a pole of
order one).
\qed

%%%%%%%%%%%%%%%

\section{Another example and a conjecture}\label{sec3}
Let $\gamma\in(0,1)$. Given a sequence 
$\{x_n\}$ with $x_{n+1}\ge x_n+h_n$, where $h_n=x_n^\gamma$, put
\[ f(x)=\left\{ \begin{array}{cc} x_i-h_i & \mathrm{if}\ x\in[x_i-h_i,x_i] \\ x_i+h_i &
    \mathrm{if}\ x\in(x_i,x_i+h_i] \\ x & \mathrm{otherwise} \end{array}\right.
\]
Again, $f$ is non-negative, non-decreasing and both $O(x^\gamma)$ and $\Omega(x^\gamma)$. 
The abscissas $\sigma_c,\sigma_a$ of convergence of $\int_1^\infty(f(x)-x)x^{-s-1}dx$ satisfy
$\sigma_c\ge 2\gamma-1$ by Proposition \ref{prop-Omega}, while comparison of $f$ with the function
considered in the proof of Proposition \ref{prop1} gives $\sigma_a\le 2\gamma-1$. This 
implies $\sigma_c=\sigma_a=2\gamma-1$. 

With $g(x)=f(x)-x$, we have 
\[ G(s)=\int_1^\infty g(x)x^{-s-1}dx = \sum_{i=1}^\infty\left(
    \int_{x_i-h_i}^{x_i} \frac{x_i-h_i-t}{t^{s+1}}dt    +\int_{x_i}^{x_i+h_i}\frac{x_i+h_i-t}{t^{s+1}}dt
   \right).\]
Since $g$ assumes positive and negative values, we must argue more carefully than above. Focusing on
a summand for fixed $i$, we have 
\bea \lefteqn{  \int_{x_i-h_i}^{x_i} \frac{x_i-h_i-t}{t^{s+1}}dt
    +\int_{x_i}^{x_i+h_i}\frac{x_i+h_i-t}{t^{s+1}}dt } \nn\\
    &&= \frac{x_i-h_i}{s}\left(\frac{1}{(x_i-h_i)^s}-\frac{1}{x_i^s}\right)
     + \frac{x_i+h_i}{s}\left(\frac{1}{x_i^s}-\frac{1}{(x_i+h_i)^s}\right)
    -\frac{1}{1-s}\left( (x_i+h_i)^{1-s}-(x_i-h_i)^{1-s}\right) \nn\\
   &&= \frac{2h_i}{sx_i^s} +\frac{1}{s(1-s)}\left( (x_i-h_i)^{1-s}-(x_i+h_i)^{1-s}\right) \nn\\
   &&= \frac{2h_i}{sx_i^s} +\frac{x_i^{1-s}}{s(1-s)}\left( (1-\frac{h_i}{x_i})^{1-s}-(1+\frac{h_i}{x_i})^{1-s}\right).\label{ft2}
\eea
In view of $h_i<<x_i$, we expand the term in the large brackets using the binomial series:
\bean  \lefteqn{ (1-\frac{h_i}{x_i})^{1-s}-(1+\frac{h_i}{x_i})^{1-s} } \\
   &&= \left( 1-(1-s)\frac{h_i}{x_i}+\frac{(1-s)(-s)}{2}(\frac{h_i}{x_i})^2-\frac{(1-s)(-s)(-s-1)}{3!}(\frac{h_i}{x_i})^3\right)\\
   && -\left( 1+(1-s) \frac{h_i}{x_i}+\frac{(1-s)(-s)}{2}(\frac{h_i}{x_i})^2+\frac{(1-s)(-s)(-s-1)}{3!}(\frac{h_i}{x_i})^3\right)
     +O((\frac{h_i}{x_i})^5)\\
   &&= -2(1-s)\frac{h_i}{x_i}-\frac{1}{3}(1-s)s(s+1)  (\frac{h_i}{x_i})^3 +O((\frac{h_i}{x_i})^5).
\eean
Plugging this into (\ref{ft2}), the first order terms cancel and we get
%\[  \int_{x_i-h_i}^{x_i} \frac{x_i-h_i-t}{t^{s+1}}dt    +\int_{x_i}^{x_i+h_i}\frac{x_i+h_i-t}{t^{s+1}}dt 
%    = -\frac{1}{3}(s+1)x_i^{1-s}\left( (\frac{h_i}{x_i})^3   +O((\frac{h_i}{x_i})^5)\right),\]
\[  \int_{x_i-h_i}^{x_i+h_i} g(x)x^{-s-1}dx= -2(s+1)x^{1-s}\left( \frac{1}{3!}(\frac{h_i}{x_i})^3+\frac{(s+2)(s+3)}{5!}(\frac{h_i}{x_i})^5+\cdots
\right).\]
With $h_i=x_i^\gamma$, where $\gamma\in(0,1)$, we have
\bean \int_1^\infty g(x)x^{-s-1}dx &=&
     -\frac{1}{3}(s+1)\sum_{i=1}^\infty x_i^{1-s}\left((\frac{x_i^\gamma}{x_i})^3 +O((\frac{h_i}{x_i})^5\right)\\
   &=& -\frac{1}{3}(s+1)\sum_{i=1}^\infty x_i^{3\gamma-s-2} +O\left(\sum_i x_i^{5\gamma-s-4}\right).
\eean
With $x_i=2^i$, the leading term equals 
$-\frac{1}{3}(s+1)\sum_{i=1}^\infty 2^{(3\gamma-s-2)i}$. From this it follows that 
the series converges if $\Re(3\gamma-s-2)<0$, or $\Re\,s>3\gamma-2$, while the sum over the higher
order terms converges for $\Re\,s>5\gamma-4$. In view of
$\sum_{i=1}^\infty 2^{(3\gamma-s-2)i}=\frac{1}{1-2^{3\gamma-s-2}}-1$,
$G$ is meromorphic on $\{\Re\,s>5\gamma-4\}$ with first order poles at 
$3\gamma-2+i\frac{2\pi}{\log 2}\7Z$.

For $\Re\,s>2\gamma-1=\sigma_a$, it is clear that the above series converges to
$\int_1^\infty(f(x)-x)x^{-s-1}dx$, so that the series gives an analytic continuation of the
integral to $\{\Re\,s>3\gamma-2\}$. We thus have:

\bprop \label{prop3}
For every $\gamma\in(0,1)$ there exists a function $f:[1,\infty)\rarr\7R$ such that 
\begin{itemize}
\item $f$ is non-negative and non-decreasing,
\item $f(x)-x$ is $O(x^\gamma)$ and $\Omega(x^\gamma)$,
\item $F(s)=\int_1^\infty f(x)x^{-s-1}dx$ converges for $s>1$,
\item $G(s)=\int_1^\infty (f(x)-x)x^{-s-1}dx$ converges if and only if $s>2\gamma-1$,
\item $G(s)=F(s)-\frac{1}{s-1}$ analytically continues to $\{\Re\,s>3\gamma-2\}$ and has a 
singularity at $3\gamma-2$.
\end{itemize}
\eprop

\brem 1. In view of $3\gamma-2<2\gamma-1$, the maximal half-plane of holomorphicity of $G$ is larger
than the half-plane of convergence of the integral defining it, reflecting the fact that the
Phragm\'en-Landau theorem  does not apply to the function $g$, which is not ultimately of one sign.

2. This shows that holomorphicity of $G$ on $\{\Re\,s>\alpha\}$ is compatible with
$f(x)-Ax=O(x^\gamma)$ being true only for $\gamma\ge\frac{\alpha+2}{3}>\frac{\alpha+1}{2}$, showing
that Theorem \ref{theor-main} is false for $A>0$ if we drop the sign condition on $f(x)-Ax$. 

3. The above function $f$ was designed in such a way as to maximize cancellations in the integral
defining $G(s)$, while being non-negative, non-decreasing and $x+\Omega(x^\gamma)$. It is hard to
see how one could construct a function with the properties in the above proposition, but with $G$
extending holomorphically to a larger half-plane.  
(One could replace the numbers $x_i+h_i$ and $x_i-h_i$ by $x_i+h'_i$ and $x_i-h_i''$,
respectively, allowing $h'_i\ne h''_i$ in the hope of achieving higher order cancellations in the
above computation. But the converse happens: While the cancellation of first order terms 
still goes through, it breaks down in second order.)

4. While the above considerations provide only some evidence for the conjecture that follows, they
prove that we cannot hope for more.
\erem

\bconj
Let $f:[1,\infty)\rarr\7R$ be non-negative and non-decreasing. Assume that 
$F(s)=\int_1^\infty f(x)x^{-s-1}dx$ converges for $s>1$ and that $F(s)-\frac{A}{s-1}$ has a
holomorphic continuation to $\{\Re\,s>\alpha\}$, where $\alpha\in(0,1)$. 
Then $f(x)=Ax+O(x^{\frac{\alpha+2}{3}+\ve})$.
\econj

%%%%%%%%%%%%%%%%%%%%%%%%%%%%


\begin{thebibliography}{99}
\bibitem{bochner} S. Bochner: Ein Satz von Landau und Ikehara. Math. Zeit. {\bf 37}, 1-9 (1933).
\bibitem{CQ} D. Choimet, H. Queffelec: {\it Twelve landmarks of twentieth-century
    analysis}. Cambridge University Press, 2015. 
\bibitem{cizek} J. \v C\'i\v zek: On the proof of the prime number theorem. 
\v Casopis P\v est. Mat. {\bf 106}, 395-401 (1981).
\bibitem{ikehara} S. Ikehara: An extension of Landau's theorem in the analytical theory of numbers. 
J. of Math. Phys. (MIT) {\bf 10}, 1-12 (1931).
\bibitem{kienast} A. Kienast: Beweis des Satzes $x^{-1}\log^q(\psi(x)-x)\rarr 0$ ohne
  \"Uberschreitung der Geraden $\sigma=1$. Math. Zeit. {\bf 43}, 113-119 (1938).
%\bibitem{kor} J. Korevaar: {\it Tauberian theory. A century of developments}. Springer, 2004.
%\bibitem{MV} H. L. Montgomery, U. M. A. Vorhauer: Changes of sign of the error term in the prime
%number theorem. Funct. et Approx. {\bf XXXV}, 235- 247 (2006).
\bibitem{RR} S. G. R\'ev\'esz, A. de Roton: Generalization of the effective Wiener-Ikehara theorem. 
Int. J. Num. Th. {\bf 9}, 2091-2128 (2013).
\bibitem{tenenbaum} G. Tenenbaum: {\it Introduction to analytic and probabilistic number theory}. Third
edition. American Mathematical Society, 2015.
\bibitem{t} G. Tenenbaum: Private communication.
\end{thebibliography}
\end{document}